\documentclass[12pt]{article}
\usepackage{amssymb, amsfonts}
\begin{document}

\title{On the reducibility of Schlesinger isomonodromic families}

\author{R.R~Gontsov and V.P.~Leksin}
\date{}

\maketitle

{\small\baselineskip 10pt  \abstract{We obtain some sufficient
conditions for reducibility of a Schlesinger isomonodromic family
with the (block) upper-triangular monodromy to the same (block)
upper-triangular form via a constant gauge transformation. We also
obtain integral representations of hypergeometric type for entries
of upper-triangular solutions of the Schlesinger equation.}}

\section{Introduction}

Let us consider on the Riemann sphere $\overline{\mathbb
C}=\mathbb C\cup \{z=\infty\}$ a \emph{Fuchsian} system of $p$
linear differential equations
If the infinity is not contained in a set of singularities of the
coefficient matrix, then such a system can be written in the form
\begin{eqnarray}\label{fuchs}
\frac{dy}{dz}=\biggl(\sum_{i=1}^n\frac{B_i^0}{z-a_i^0}\biggr)y,
\qquad \sum_{i=1}^nB_i^0=0,
\end{eqnarray}
where $y(z)\in{\mathbb C}^p,$ $B_1^0,\ldots,B_n^0$ are constant
matrices of size $p\times p$ (they are residue matrices of the
system), $a_1^0,\ldots,a_n^0\in{\mathbb C}$ are singular points of
the system.

An important characteristic of a linear system is the monodromy,
or monodromy representation. In a neighborhood of a non-singular
point $z_0$ consider a fundamental matrix $Y(z)$ of the system
(\ref{fuchs}). An analytic continuation of the matrix $Y(z)$ along
an arbitrary loop $\gamma$ starting at $z_0$ and contained in
$\overline{\mathbb C}\setminus \{a_1^0,\ldots,a_n^0\}$ transforms
this matrix in, generally speaking, another fundamental matrix
$\widetilde Y(z)$. Two fundamental matrices are connected by a
non-degenerate transition matrix $G_{\gamma}$ corresponding to the
loop $\gamma$:
$$
\widetilde Y(z)=Y(z)G_{\gamma}.
$$
The map $[\gamma]\mapsto G_{\gamma}^{-1}$ depends on the homotopic
class  $[\gamma]$ of the loop $\gamma$ only and thus defines a
representation
$$
\chi: \pi_1(\overline{\mathbb C}\setminus\{a_1^0,\ldots,a_n^0\},
z_0)\longrightarrow{\rm GL}(p,{\mathbb C})
$$
of the fundamental group of the space $\overline{\mathbb
C}\setminus \{a_1^0,\ldots,a_n^0\}$ in the space of non-degenerate
complex matrices of size  $p\times p$. The representation $\chi$
is called the {\it monodromy} of the system (\ref{fuchs}).

By the {\it monodromy matrix} of the Fuchsian system (\ref{fuchs})
at a singular point $a_i^0$ (with respect to the fundamental
matrix $Y(z)$) one understands the matrix $G_i$ corresponding to a
simple loop $\gamma_i$ encircling the point $a_i^0$, so that
$G_i^{-1}=\chi([\gamma_i])$. The matrices $G_1,\ldots,G_n$ are
generators of the {\it monodromy group} of the system
(\ref{fuchs}).

If instead of the fundamental matrix $Y(z)$ we consider another
fundamental matrix $Y'(z)=Y(z)C$, $C\in{\rm GL}(p, {\mathbb C})$,
then the corresponding monodromy matrices are of the form
$G'_i=C^{-1}G_iC$. Thus monodromy matrices are defined up to a
simultaneous conjugation.

The main object of the present work is a {\it Schlesinger
isomonodromic family}
\begin{eqnarray}\label{schl}
\frac{dy}{dz}=\biggl(\sum_{i=1}^n\frac{B_i(a)}{z-a_i}\biggr)y,
\qquad B_i(a^0)=B_i^0,
\end{eqnarray}
of Fuchsian systems holomorphically depending on a parameter
$a=(a_1,\ldots,a_n)\in D(a^0)$, where $D(a^0)$ is a small polydisk
centered at the point $a^0=(a_1^0,\ldots,a_n^0)$ of the space
${\mathbb C}^n\setminus\bigcup_{i\ne j}\{a_i=a_j\}$. The monodromy
of systems of this family is the same for all $a\in D(a^0)$ and
the residue matrices $B_i(a)$ satisfy the {\it Schlesinger
equation} \cite{Sch}
$$
dB_i(a)=-\sum_{j=1, j\ne i}^n\frac{[B_i(a),B_j(a)]}{a_i-a_j}\,
d(a_i-a_j), \qquad i=1,\ldots,n.
$$

The Schlesinger equation is {\it integrable in the Frobenius
sense} in the polydisk $D(a^0)$, i.e., for any initial data
$B_1^0,\ldots,B_n^0$ it has the unique solution
$B_1(a),\ldots,B_n(a)$ such that $B_i(a^0)=B_i^0$ (see
\cite[theorem 14.2]{Bo}). According to the Malgrange theorem
\cite{Ma} the matrix valued functions $B_i(a)$ can be continued
meromorphically onto the universal cover $Z$ of the space
${\mathbb C}^n\setminus\bigcup_{i\ne j}\{a_i=a_j\}$. The polar set
$\Theta\subset Z$ of the continued matrix functions $B_i(a)$ is
called the {\it Malgrange $\Theta$-divisor} ($\Theta$ depends on
the initial data $B_i(a^0)=B_i^0$). Moreover there exists a
function $\tau$, holomorphic on the whole space $Z$, whose zero
set coincides with $\Theta$. This function is called the {\it
$\tau$-function} of the Schlesinger equation and Miwa's formula
\cite{JMU} (see also \cite[theorem 17.1]{Bo}) holds for it:
$$
d\ln\tau(a)=\frac12\sum_{i=1}^n\sum_{j=1, j\ne i}^n\frac{{\rm tr}
(B_i(a)B_j(a))}{a_i-a_j}\,d(a_i-a_j).
$$

If the system (\ref{fuchs}) is meromorphically equivalent to a
Fuchsian system
$$
\frac{d\tilde y}{dz}=\biggl(\sum_{i=1}^n\frac{\widetilde B_i}
{z-a_i^0}\biggr)\tilde y, \qquad \sum_{i=1}^n\widetilde B_i=0
$$
(under a meromorphically invertible gauge transformation $\tilde
y=\Gamma(z)y$) whose residue matrices $\widetilde B_i$ have a
block upper-triangular form
$$
\widetilde B_i=\left(\begin{array}{cccc} \widetilde B_i^1 & & * & * \\
                                       0 & \widetilde B_i^2 & & * \\
                                       \vdots & \ddots & \ddots & \\
                                       0 & \ldots & 0 & \widetilde B_i^k
                     \end{array}\right), \qquad i=1,\ldots,n,
$$
then the monodromy matrices of the (transformed and, consequently,
initial) system have the same block upper-triangular form. The
inverse, generally speaking, is fulfilled in a weaker sense.
Namely, as shown by S.~Malek \cite{Ma1}, if the monodromy matrices
of the Fuchsian system (\ref{fuchs}) have a block upper-triangular
form
\begin{eqnarray}\label{blocktr}
G_i=\left(\begin{array}{cccc} G_i^1 & & * & * \\
                              0 & G_i^2 & & * \\
                         \vdots & \ddots & \ddots & \\
                              0 & \ldots & 0 & G_i^k
                     \end{array}\right), \qquad i=1,\ldots,n,
\end{eqnarray}
then the system is meromorphically equivalent to a Fuchsian system
whose residue matrices $\widetilde B_i$ are of the block
upper-triangular form
$$
\widetilde B_i=\left(\begin{array}{cc} B'_i & * \\
                                       0 & B''_i \\
                     \end{array}\right).
$$
But a subsequent simultaneous reducibility of the blocks $B'_i$
(or $B''_i$) in a correspondence with the reducibility
(\ref{blocktr}) of the matrices $G_i$ can not to take place (see
\cite[proposition 5.1.1]{Bo2}, \cite{G2}).

Similar connections between the meromorphic reducibility of
systems and the block structure of monodromy matrices take place
for Schlesinger isomonodromic families as well (more details on
the parametric reducibility will be presented elsewhere). Here,
for subsequent applications to integral representations of
solutions of the Schlesinger equation, we are interested in
sufficient conditions for the reducibility of the residue matrices
$B_i(a)$ of the Schlesinger isomonodromic family (\ref{schl}) to
the same block upper-triangular form (\ref{blocktr}) as its
monodromy matrices do have, with respect to a {\it constant} gauge
transformation. As a technique tool we employ the approach that
uses holomorphic vector bundles and meromorphic connections.

\section{Fuchsian systems and logarithmic connections in holomorphic vector bundles}

According to the Levelt theorem  ~\cite{Le}, in a neighborhood of
each singular point $a_i^0$ of the system (\ref{fuchs}) there
exists a fundamental matrix of the form
\begin{eqnarray}\label{levelt}
Y(z)=U_i(z)(z-a_i^0)^{\Lambda_i}(z-a_i^0)^{E_i},
\end{eqnarray}
where $U_i(z)$ is a holomorphically invertible matrix at the point
$a_i^0$, $\Lambda_i={\rm diag}(\lambda_i^1,\ldots,\lambda_i^p)$ is
a diagonal integer matrix whose elements $\lambda_i^j$ organize a
non-increasing sequence, $E_i=(1/2\pi{\bf i})\ln G_i$ is an
upper-triangular matrix (the normalized logarithm of the
corresponding monodromy matrix) whose eigenvalues $\rho_i^j$
satisfy the condition
$$
0\leqslant{\rm Re}\,\rho_i^j<1.
$$
Such the fundamental matrix is called the {\it Levelt matrix}, and
one also says that its columns form the {\it Levelt basis} in the
solution space of the Fuchsian system (in a neighborhood of the
singular point $a_i^0$). The diagonal elements $\lambda_i^j$ of
the matrix $\Lambda_i$ are called the (Levelt) {\it valuations}
and the complex numbers $\beta_i^j=\lambda_i^j+\rho_i^j$ are
called the (Levelt) {\it exponents} of the Fuchsian system at the
singular point $a_i^0$. It is not difficult to check that the
exponents of the Fuchsian system at the point $a_i^0$ coincide
with the eigenvalues of the residue matrix $B_i^0$.

Let us recall some notions concerning holomorphic vector bundles
and meromorphic connections. In an analytic interpretation, a
holomorphic bundle $E$ of rank $p$ (over the Riemann sphere) is
defined by {\it a cocycle} $\{g_{\alpha\beta}(z)\}$, that is, a
collection of holomorphic matrix functions, corresponding to a
covering $\{U_{\alpha}\}$ of the Riemann sphere:
$$
g_{\alpha\beta}: U_{\alpha}\cap U_{\beta}\longrightarrow{\rm GL}
(p,\mathbb C), \qquad U_{\alpha}\cap U_{\beta}\ne\varnothing.
$$
These functions satisfy conditions
$$
g_{\alpha\beta}=g_{\beta\alpha}^{-1}, \qquad
g_{\alpha\beta}g_{\beta\gamma}g_{\gamma\alpha}=I \quad (\mbox{for
} U_{\alpha}\cap U_{\beta}\cap U_{\gamma}\ne\varnothing).
$$
Two holomorphically equivalent cocycles $\{g_{\alpha\beta}(z)\}$,
$\{g'_{\alpha\beta}(z)\}$ define the same bundle. The equivalence
of cocycles means that there exists a set $\{h_{\alpha}(z)\}$ of
holomorphic functions $h_{\alpha}: U_{\alpha}\longrightarrow{\rm
GL}(p, {\mathbb C})$ such that
\begin{eqnarray}\label{equiv}
h_{\alpha}(z)g_{\alpha\beta}(z)=g'_{\alpha\beta}(z)h_{\beta}(z).
\end{eqnarray}
A {\it section} $s$ of the bundle $E$ is determined by a set
$\{s_{\alpha}(z)\}$ of vector functions $s_{\alpha}:
U_{\alpha}\longrightarrow{\mathbb C}^p$ that satisfy the
conditions $s_{\alpha}(z)=g_{\alpha\beta}(z) s_{\beta}(z)$ on
intersections $U_{\alpha}\cap U_{\beta}\ne\varnothing$ .

A {\it meromorphic connection} $\nabla$ in a holomorphic vector
bundle $E$ is determined by a set $\{\omega_{\alpha}\}$ of matrix
meromorphic differential 1-forms that are defined in corresponding
neighborhoods $U_{\alpha}$ and satisfy gluing conditions
\begin{eqnarray}\label{glue}
\omega_{\alpha}=(dg_{\alpha\beta})g_{\alpha\beta}^{-1}+
g_{\alpha\beta}\omega_{\beta}g_{\alpha\beta}^{-1}\qquad (\mbox{for
} U_{\alpha}\cap U_{\beta}\ne\varnothing).
\end{eqnarray}
Under a transition to an equivalent cocycle $\{g'_{\alpha\beta}\}$
connected with the initial one by the relations (\ref{equiv}), the
1-forms $\omega_{\alpha}$ of the connection  $\nabla$ are
transformed into the corresponding 1-forms
\begin{eqnarray}\label{formequiv}
\omega'_{\alpha}=(dh_{\alpha})h_{\alpha}^{-1}+
h_{\alpha}\omega_{\alpha}h_{\alpha}^{-1}.
\end{eqnarray}
Inversely, the existence of holomorphic matrix functions
$h_{\alpha}: U_{\alpha}\longrightarrow{\rm GL}(p, {\mathbb C})$
such that the matrix 1-forms $\omega_{\alpha}$ and
$\omega'_{\alpha}$ (satisfying conditions (\ref{glue}) for
$g_{\alpha\beta}$ and $g'_{\alpha\beta}$ respectively) are
connected by the relation (\ref{formequiv}), indicates the
equivalence of cocycles $\{g_{\alpha\beta}\}$ and
$\{g'_{\alpha\beta}\}$.

A set $\{s_{\alpha}(z)\}$ of vector functions satisfying linear
differential equations  $ds_{\alpha}=\omega_{\alpha} s_{\alpha}$
in the corresponding $U_{\alpha}$, by virtue of conditions
(\ref{glue}) determines a section of the bundle $E$, which is
called {\it horizontal} with respect to the connection $\nabla$.
Thus horizontal sections of a holomorphic vector bundle with a
meromorphic connection are determined by solutions of local
systems of linear differential equations. The {\it monodromy of a
connection} characterizes ramification of horizontal sections
under their analytic continuation along loops in
$\overline{\mathbb C}$ not containing singular points of the
connection 1-form. It is defined similarly to the monodromy of a
linear differential system. A connection is called {\it
logarithmic} or {\it Fuchsian}, if all singular points of its
1-form are poles of the first order.

The degree $\deg E$ (which is an integer) of the holomorphic
vector bundle $E$ with the meromorphic connection $\nabla$ one can
define as the sum
$$
\deg E=\sum_{i=1}^n{\rm res}_{a_i^0}\,{\rm tr}\,\omega_i
$$
of residues of local 1-forms ${\rm tr}\,\omega_i$ by the all
singular points of the connection, where $\omega_i$ is the local
1-form of the connection $\nabla$ in a neighborhood of its
singular point $a_i^0$.

Further we say about holomorphic vector bundles with logarithmic
connections in view of their close relation with Fuchsian systems.
If a bundle is holomorphically trivial (all matrices of the
cocycle can be taken as identity matrices), then by virtue of the
conditions (\ref{glue}) the matrix 1-forms of a logarithmic
connection coincide on non-empty intersections $U_{\alpha}\cap
U_{\beta}$. Hence horizontal sections of such a bundle are
solutions of a global Fuchsian system of linear differential
equations defined on the whole Riemann sphere. Inversely, the
Fuchsian system (\ref{fuchs}) determines the logarithmic
connection in the holomorphically trivial vector bundle of rank
$p$ over $\overline{\mathbb C}$. Clear, such bundle has the
standard definition by the cocycle that consists of the identity
matrices while the connection is defined by the matrix 1-form
$$
\omega^0=\sum_{i=1}^n\frac{B_i^0}{z-a_i^0}\,dz
$$
of coefficients of the system. But for us it will be more
convenient to use the following coordinate description.

At first we consider a covering  $\{U_{\alpha}\}$ of the punctured
Riemann sphere $\overline{\mathbb C}\setminus
\{a_1^0,\ldots,a_n^0\}$ and a corresponding set of constant matrix
functions $g'_{\alpha\beta}(z)\equiv{\rm const}$, which are
expressed in terms of the monodromy matrices $G_1,\ldots,G_n$ of
the system (\ref{fuchs}) via operations of multiplication and
taking the inverse matrix (see \cite[Lect. 8]{Bo}). In this case
matrix differential 1-forms $\omega'_{\alpha}$ defining a
connection will be equal to zero. Further the covering
$\{U_{\alpha}\}$ is complemented by small neighborhoods $O_i$ of
the singular points $a_i^0$ of the system, thus we obtain a
covering of the Riemann sphere $\overline{\mathbb C}$. To
non-empty intersections $O_i\cap U_{\alpha}$ there correspond
matrix functions $g'_{i\alpha}(z)=Y_i(z)$ of the cocycle, where
$Y_i(z)$ is a germ of a fundamental matrix of the system whose
monodromy matrix at the point $a_i^0$ is equal to $G_i$. (So, for
the analytic continuations of a chosen germ to non-empty
intersections $O_i\cap U_{\alpha}\cap U_{\beta}$ the cocycle
relations $g_{i\alpha}g_{\alpha\beta}=g_{i\beta}$ hold.) Matrix
differential 1-forms $\omega'_i$ determining the connection in
neighborhoods $O_i$ coincide with the 1-form $\omega^0$ of
coefficients of the system. In order to prove holomorphic
equivalence of the cocycle $\{g'_{\alpha\beta}, g'_{i\alpha}\}$ to
the identity cocycle, it remains to check existence of holomorphic
matrix functions
$$
h_{\alpha}: U_{\alpha}\longrightarrow{\rm GL}(p, {\mathbb C}),
\qquad h_i: O_i\longrightarrow{\rm GL}(p, {\mathbb C}),
$$
such that
\begin{eqnarray}\label{formequiv2}
\omega'_{\alpha}=(dh_{\alpha})h_{\alpha}^{-1}+
h_{\alpha}\omega_{\alpha}h_{\alpha}^{-1}, \qquad
\omega'_i=(dh_i)h_i^{-1}+ h_i\omega_ih_i^{-1}.
\end{eqnarray}

Since we have $\omega_{\alpha}=\omega^0$ and $\omega'_{\alpha}=0$
for all $\alpha$, the first equation of the (\ref{formequiv2}) is
rewritten as a linear system
$$
d(h^{-1}_{\alpha})=\omega^0h^{-1}_{\alpha},
$$
which has a holomorphic solution $h^{-1}_{\alpha}: U_{\alpha}
\longrightarrow{\rm GL}(p,{\mathbb C})$ because the 1-form
$\omega^0$ is holomorphic in $U_{\alpha}$. The second equation of
the (\ref{formequiv2}) has a holomorphic solution $h_i(z)\equiv
I$, as $\omega_i=\omega'_i=\omega^0$.
\medskip

One says that the bundle $E$ has a subbundle $E'\subset E$ of rank
$k<p$ that is {\it stabilized} by the connection $\nabla$, if the
pair $(E, \nabla)$ admits a coordinate description
$\{g_{\alpha\beta}\}$, $\{\omega_{\alpha}\}$ of the following
block-uppertriangular form:
$$
g_{\alpha\beta}= \left(\begin{array}{cc}
       g_{\alpha\beta}^1 & * \\
       0 & g_{\alpha\beta}^2 \\
      \end{array}\right), \qquad
\omega_{\alpha}= \left(\begin{array}{cc}
       \omega_{\alpha}^1 & * \\
       0 & \omega_{\alpha}^2 \\
      \end{array}\right),
$$
where $g_{\alpha\beta}^1$ and $\omega_{\alpha}^1$ are blocks of
size $k\times k$ (then the cocycle $\{g_{\alpha\beta}^1\}$ defines
the subbundle $E'$ and the 1-forms $\omega_{\alpha}^1$ define the
restriction $\nabla'$ of the connection $\nabla$ to the subbundle
$E'$).
\medskip

{\bf Example 1.} Consider the Fuchsian system (\ref{fuchs}) with a
reducible monodromy representation and corresponding
holomorphically trivial vector bundle with the logarithmic
connection. Let us demonstrate that to the monodromy
subrepresentation there corresponds a holomorphic vector subbundle
that is stabilized by the connection.

We use the above coordinate description of the bundle and
connection with the cocycle $\{g'_{\alpha\beta}, g'_{i\alpha}\}$
and set $\{\omega'_{\alpha}, \omega'_i\}$ of matrix 1-forms. We
can pass to the equivalent cocycle (preserving the previous
notations) changing the matrices $g'_{\alpha\beta}$ to the
matrices $S^{-1}g'_{\alpha\beta}S$ and the matrices $g'_{i\alpha}$
to the matrices $g'_{i\alpha}S$, where $S$ is a constant
non-degenerate matrix reducing the monodromy matrices
$G_1,\ldots,G_n$ of the system to the same block upper-triangular
form $G'_1,\ldots,G'_n$. Then the matrices $g'_{\alpha\beta}$ are
block upper-triangular (and $\omega'_{\alpha}=0$), and the
matrices $g'_{i\alpha}$ have the form
$$
g'_{i\alpha}(z)=M_i(z)(z-a_i^0)^{E_i}, \qquad E_i=(1/2\pi{\bf
i})\ln G'_i,
$$
where $M_i(z)$ are meromorphic matrices in neighborhoods of the
corresponding points $a_i^0$. For the latter the following
factorizations hold: $M_i(z)=V_i(z)P_i(z)$, where the matrix
$V_i(z)$ is holomorphically invertible at the point $a_i^0$,
$P_i(z)$ is a polynomial {\it upper-triangular} matrix in
$(z-a_i^0)^{\pm1}$ (see, for example, ~\cite[Lemma 1]{Go}). Thus
changing the matrices $g'_{i\alpha}(z)$ to
$$
V^{-1}_i(z)g'_{i\alpha}(z)=P_i(z)(z-a_i^0)^{E_i}
$$
and the matrix 1-forms $\omega'_i$ to
$$
V^{-1}_i\omega'_iV_i-V^{-1}_i(dV_i)=(dP_i)P^{-1}_i
+P_i\frac{E_idz}{z-a_i^0}P^{-1}_i,
$$
we pass to the holomorphically equivalent coordinate description
whose cocycle matrices and matrix 1-forms of the connection have
the same block upper-triangular form.
\medskip

{\bf Remark 1.} From the definition of the subbundle $E'$ that is
stabilized by the logarithmic connection $\nabla$, one can see
that the set of exponents of the restriction $\nabla'$ at each
singular point $a_i^0$ (the set of eigenvalues of the residue
matrix of the corresponding 1-form) is a subset of exponents of
the connection $\nabla$. Therefore the degree $\deg E'$ of the
subbundle $E'$ is equal to the sum of exponents from these subsets
over all singular points of the connection.
\medskip

The following auxiliary lemma points to a certain block structure
of the residue matrices of a Fuchsian system in the case when the
corresponding holomorphically trivial vector bundle with a
logarithmic connection has a {\it holomorphically trivial}
subbundle that is stabilized by the connection.
\medskip

{\bf Lemma 1.} {\it If $E$ is the holomorphically trivial vector
bundle of rank $p$ over $\overline{\mathbb C}$ having a
holomorphically trivial subbundle $E'\subset E$ of rank $k$ that
is stabilized by the connection $\nabla$, then the corresponding
Fuchsian system $(\ref{fuchs})$ is reduced to a block
upper-triangular form via a constant gauge transformation $\tilde
y(z)=Cy(z)$, $C\in{\rm GL} (p,\mathbb C)$. That is,
$$
CB_i^0C^{-1}=\left(\begin{array}{cc} B'_i & * \\
                                        0 & *
                   \end{array}\right),\qquad i=1,\ldots,n,
$$
where $B'_i$ is a block of size $k\times k$.}
\medskip

{\bf Proof.} Let $\{s_1,\ldots,s_p\}$ be a basis of global
sections of the bundle $E$ (which are linear independent at each
point $z\in\overline{\mathbb C}$) such that the 1-form of the
connection $\nabla$ in this basis is the 1-form $\omega^0$ of
coefficients of the Fuchsian system. Consider also a basis
$\{s'_1,\ldots,s'_p\}$ of global holomorphic sections of the
bundle $E$ such that $s'_1,\ldots,s'_k$ are sections of the
subbundle $E'$, $(s'_1,\ldots,s'_p)=(s_1,\ldots,s_p)C^{-1}$,
$C\in{\rm GL}(p,{\mathbb C})$.

Now choose a basis $\{h_1,\ldots,h_p\}$ of sections of the bundle
$E$ such that they are horizontal with respect to the connection
$\nabla$ and $h_1,\ldots,h_k$ are sections of the subbundle $E'$
(it is possible since $E'$ is stabilized by the connection). Let
$Y(z)$ be a fundamental matrix of the Fuchsian system whose
columns are the coordinates of the sections $h_1,\ldots,h_p$ in
the basis $\{s_1,\ldots,s_p\}$. Then
$$
\widetilde Y(z)=CY(z)=\left(\begin{array}{cc} k\times k & * \\
                                            0 & *
                      \end{array}\right)
$$
is a block upper-triangular matrix, since its columns are the
coordinates of the sections $h_1,\ldots,h_p$ in the basis
$\{s'_1,\ldots, s'_p\}$. Consequently, the transformation $\tilde
y(z)=Cy(z)$ reduces the initial system to a block upper-triangular
form. {\hfill $\Box$}

\section{On the reducibility of Fuchsian systems and their isomonodromic families}

Let monodromy matrices of the Fuchsian system (\ref{fuchs}) have
the block upper-triangular form (\ref{blocktr}). Denote by
$m^s\times m^s$ size of the blocks $G^s_1,\ldots,G^s_n$
($s=1,\ldots,k$). There holds the following sufficient condition
of reducibility of the residue matrices of the system to the same
block upper-triangular form.
\medskip

{\bf Theorem 1.} {\it If the exponents $\beta_i^j$ of the Fuchsian
system $(\ref{fuchs})$ satisfy the condition
\begin{eqnarray}\label{ineq1}
{\rm Re}\,\beta_i^j>-1/n(p-m^k), \qquad i=1,\ldots,n, \quad
j=1,\ldots,p,
\end{eqnarray}
then there exists a constant matrix $C\in{\rm GL}(p,\mathbb C)$
such that the matrices $CB_i^0C^{-1}$ have the same block
upper-triangular form as the monodromy matrices of the system.}
\medskip

{\bf Proof.} We use a geometric interpretation (exposed in the
previous section) according to which to the Fuchsian system
(\ref{fuchs}) there corresponds a holomorphically trivial vector
bundle $E$ of rank $p$ over the Riemann sphere endowed a
logarithmic connection $\nabla$. Since the monodromy matrices of
the system are block upper-triangular, there exists a flag
$E^1\subset E^2\subset\ldots\subset E^k=E$ of subbundles of ranks
$m^1,m^1+m^2,\ldots,m^1+\ldots+m^k=p$ correspondingly that are
stabilized by the connection $\nabla$ (see Example 1).

Let us estimate the degree of each subbundle $E^s$, $s\leqslant
k$, using Remark 1. The degree of the holomorphically trivial
vector bundle $E^k$ is equal to zero and for $s<k$ we have:
\begin{eqnarray*}
\deg E^s = \sum_{i=1}^n\sum_{j\in J_i,\,|J_i|=m^1+\ldots+m^s}{\rm
Re}\,\beta_i^j>-(m^1+\ldots+m^s)/(p-m^k)\geqslant-1.
\end{eqnarray*}
Therefore, $\deg E^s=0$ (the degree of a subbundle of a
holomorphically trivial vector bundle is non-positive, see
\cite[Prop. 11.1]{Bo}), and all the subbundles
$E^1\subset\ldots\subset E^k$ are holomorphically trivial (a
subbundle of a holomorphically trivial vector bundle is
holomorphically trivial, if its degree is equal to zero, see
\cite[Corollary 11.1]{Bo}). Now the assertion of the theorem
follows from Lemma 1. {\hfill $\Box$}
\medskip

{\bf Remark 2.} Although the inequalities (\ref{ineq1}) point out
to the boundedness of real parts of the exponents from below, from
these inequalities and the classical Fuchs relation \cite{Le}
$$
\sum_{i=1}^n\sum_{j=1}^p\beta_i^j=0
$$
(which says that the degree of the holomorphically trivial vector
bundle $E$ is equal to zero) it follows also the boundedness from
above.
\medskip

{\bf Corollary 1.} {\it If the monodromy representation of the
Fuchsian system $(\ref{fuchs})$ is upper-triangular and its
exponents $\beta_i^j$ satisfy the condition
\begin{eqnarray*}
{\rm Re}\,\beta_i^j>-1/n(p-1), \qquad i=1,\ldots,n, \quad
j=1,\ldots,p,
\end{eqnarray*}
then there exists a constant matrix $C\in{\rm GL}(p,\mathbb C)$
such that all the matrices $CB_i^0C^{-1}$ are upper-triangular.}
\medskip

The monodromy representation of a Fuchsian system is called a {\it
B-representation}, if it is reducible and the Jordan form of each
monodromy matrix $G_i$ consists of one Jordan box only. In this
case we have the following assertion.
\medskip

{\bf Proposition 1.} {\it If the monodromy representation of the
Fuchsian system $(\ref{fuchs})$ is an upper-triangular
B-representation, then there exists a constant matrix $C\in{\rm
GL}(p,\mathbb C)$ such that all the matrices $CB_i^0C^{-1}$ are
upper-triangular.}
\medskip

{\bf Proof.} Since the monodromy representation of the Fuchsian
system is a B-representation, then at each singular point $a_i$
the system has only one exponent (of multiplicity $p$):
$\beta_i^1=\ldots=\beta_i^p=\beta_i$ (see the proof of Theorem
11.2 from \cite{Bo}). Then we have $p\sum_{i=1}^n\beta_i=0$
because of triviality of the vector bundle $E$ of rank $p$
corresponding to the given Fuchsian system. Now, as in the proof
of Theorem 1, from the upper-triangularity of the monodromy
matrices it follows the existence of a flag $E^1\subset E^2\subset
\ldots\subset E^p=E$ of subbundles of ranks $1, 2, \ldots, p$
correspondingly that are stabilized by the connection $\nabla$.
But each subbundle $E^s$ is holomorphically trivial, since its
degree $\deg E^s= s\sum_{i=1}^n\beta_i$ is equal to zero. To
complete the proof it remains to apply again Lemma 1. {\hfill
$\Box$}
\medskip

Let us say a few words about the reducibility of Schlesinger
isomonodromic families (\ref{schl}). First of all, it follows from
the result of S.\,Malek (mentioned in Introduction) that if the
monodromy of such a family is reducible, then via a
meromorphically invertible (in $z$) gauge transformation $\tilde
y=\Gamma(z,a)y$ this family can be transformed to a Schlesinger
isomonodromic family whose residue matrices $\widetilde B_i(a)$
have the block-uppertriangular form
$$
\widetilde B_i(a)=\left(\begin{array}{cc} B'_i(a) & * \\
                                          0 & B''_i(a) \\
                     \end{array}\right).
$$
Indeed, for $a=a^0$ the initial system is reduced to the required
form via a meromorphically invertible transformation. The
transformed block upper-triangular system can be included in the
Schlesinger isomonodromic family with the residue matrices
$\widetilde B_i(a)$ of the same block upper-triangular form
(because they are solutions of the Schlesinger equation and they
are block upper-triangular for $a=a^0$). Since the initial and
obtained families have the same monodromy, they are connected by a
meromorphically invertible transformation $\tilde y=\Gamma(z,a)y$.
S.~Malek \cite{Ma2} has shown also that the matrix $\Gamma(z,a)$
is holomorphic with respect to the variable $a$ in $D(a^0)$, with
the exception of some analytic subset of codimension one, and the
entries of $\Gamma(z,a)$ are rational functions in $z$,
$a_1,\ldots,a_n$ and entries of residue matrices $B_i(a)$ (see
also \cite{DM}).

We can formulate the following analogs of Corollary 1 and
Proposition 1.
\medskip

{\bf Proposition 2.} {\it If the monodromy representation of the
Fuchsian system $(\ref{fuchs})$ is upper-triangular and its
exponents $\beta_i^j$ satisfy the condition
\begin{eqnarray*}
{\rm Re}\,\beta_i^j>-1/n(p-1), \qquad i=1,\ldots,n, \quad
j=1,\ldots,p,
\end{eqnarray*}
then for the Schlesinger isomonodromic deformation $(\ref{schl})$
of this system there exists a constant matrix $C\in{\rm GL}
(p,\mathbb C)$ such that all the matrices $CB_i(a)C^{-1}$ are
upper-triangular.}
\medskip

{\bf Proposition 3.} {\it If the monodromy representation of the
Fuchsian system $(\ref{fuchs})$ is an upper-triangular
B-representation, then for the Schlesinger isomonodromic
deformation $(\ref{schl})$ of this system there exists a constant
matrix $C\in{\rm GL}(p,\mathbb C)$ such that all the matrices
$CB_i(a)C^{-1}$ are upper-triangular.}
\medskip

{\bf Proof.} Both propositions are proved with the same method.
The residue matrices $B_i^0$ of the Fuchsian system satisfying the
condition of the proposition is transformed to an upper-triangular
form $CB_i^0C^{-1}$. It is not difficult to see that matrices
$CB_i(a)C^{-1}$ also satisfy the Schlesinger equation. As they are
upper-triangular for $a=a^0$: $CB_i(a^0)C^{-1}=CB_i^0C^{-1}$, then
by virtue of the Schlesinger equation they remain upper-triangular
for all $a\in D(a^0)$. Indeed, any partial derivative of solutions
of the Schlesinger equation is expressed via partial derivatives
of lower orders in terms of matrix operations of addition and
multiplication. Therefore from the upper-triangularity of
solutions for $a=a^0$ it follows that their partial derivatives of
any order are upper-triangular for $a=a^0$. Hence such solutions
are upper-triangular for all $a\in D(a^0)$. {\hfill $\Box$}
\medskip

Since the eigenvalues $\beta_i^j$ of the residue matrices $B_i(a)$
of the Schlesinger isomonodromic family (\ref{schl}) do not depend
on $a$ (see, for example, \cite{Bo3}), for the $\tau$-function of
an upper-triangular family (as a consequence of Miwa's formula)
there holds the relation
$$
d\ln\tau(a)=\frac12\sum_{i=1}^n\sum_{j=1, j\ne i}^n
\frac{\alpha_{ij}}{a_i-a_j}\,d(a_i-a_j),
$$
where $\alpha_{ij}=\beta_i^1\beta_j^1+\ldots+\beta_i^p\beta_j^p$.
Thus $\tau(a)=\prod_{i<j}(a_i-a_j)^{\alpha_{ij}}$ is a non-zero
holomorphic function on the universal cover $Z$ of the space
${\mathbb C}^n\setminus\bigcup_{i\ne j}\{a_i=a_j\}$. This implies
that the Malgrange $\Theta$-divisor of an {\it upper-triangular}
Schlesinger isomonodromic family is empty and one can {\it
holomorphically} continue all matrices $B_i(a)$ to the whole space
$Z$.
\medskip

\section{Triangular Schlesinger families in small dimensions}

In the previous section there were proposed some sufficient
conditions of the reducibility of the Schlesinger isomonodromic
family (\ref{schl}) to an upper-triangular form via a constant
gauge transformation. Now we describe integral representations for
entries of the residue matrices of such a family in the dimensions
$p=2$ and $p=3$.

In the case $p=2$ upper-triangular residue matrices $B_i(a)$ have
the form
$$
B_i(a)=\left(\begin{array}{cc} \beta_i^1 & b_i(a) \\
                               0 & \beta_i^2
             \end{array}\right), \qquad i=1,\ldots,n,
$$
where the exponents $\beta_i^1, \beta_i^2$ are constant. Then the
functions $b_i(a)$ by virtue of the Schlesinger equation satisfy
the following system of homogeneous linear differential equations:
\begin{eqnarray}\label{schl2}
db_i(a)=-\sum_{j=1,j\ne i}^n(\beta_ib_j(a)-\beta_jb_i(a))\,
\frac{d(a_i-a_j)}{a_i-a_j}, \qquad i=1,\ldots,n,
\end{eqnarray}
where $\beta_i=\beta_i^1-\beta_i^2$.

The system (\ref{schl2}) is the Jordan-Pochhammer system, i. e., a
linear Pfaffian system for the vector function
$b(a)=(b_1(a),\ldots,b_n(a))^{\top}$:
\begin{eqnarray}\label{JPsyst}
db=\Omega\, b, \qquad b(a)\in{\mathbb C}^n,
\end{eqnarray}
with the meromorphic (holomorphic in $D(a^0)$) coefficient matrix
1-form
\begin{eqnarray}\label{JPform}
\Omega=\sum_{1\leqslant j<k\leqslant n}J_{jk}(\beta)\,
\frac{d(a_j-a_k)}{a_j-a_k},
\end{eqnarray}
where $J_{jk}(\beta)$ are constant $(n\times n)$-matrices. Each
matrix $J_{jk}(\beta)$ has only four non-zero entries: in the
$j$-th row the entry with the number $j$ is equal to $\beta_k$
while the entry with the number $k$ is equal to $-\beta_j$, and in
the $k$-th row the entry with the number $j$ is equal to
$-\beta_k$ while the entry with the number $k$ is equal to
$\beta_j$.

Any solution of the system (\ref{JPsyst}), (\ref{JPform}) can be
represented in the following form (\cite{Ko}, \cite{KM}):
\begin{eqnarray}\label{JP}
b_i(a)=\beta_i\int_{\gamma}(t-a_1)^{\beta_1}\ldots(t-a_n)^{\beta_n}\frac{dt}{t-a_i},
\qquad i=1,\ldots,n,
\end{eqnarray}
where $\gamma\in H_1(\overline{\mathbb C}\setminus
\{a_1,\ldots,a_n\}, {\mathbb L}^*)$ is a twisted cycle whose
coefficient takes value in the local system ${\mathbb L}^*$ which
is dual to the local system $\mathbb L$ associated to the
multi-valued function $F(t)=(t-a_1)^{\beta_1}\ldots
(t-a_n)^{\beta_n}$. (To define correctly integration of
multi-valued functions along cycles, one should use 1-homologies
with coefficients in a local system; see details, for example, in
\cite{DMo} or \cite{OT}.)

In the case $p=3$ upper-triangular residue matrices $B_i(a)$ have
the form
$$
B_i(a)=\left(\begin{array}{ccc} \beta_i^1 & u_i(a) & b_i(a) \\
                                0 & \beta_i^2 & v_i(a) \\
                                0 & 0 & \beta_i^3
             \end{array}\right), \qquad i=1,\ldots,n,
$$
where the exponents $\beta_i^1, \beta_i^2, \beta_i^3$ are
constant. By virtue of the Schlesinger equation, the vector
functions $u(a)=(u_1(a),\ldots,u_n(a))^{\top}$ and
$v(a)=(v_1(a),\ldots,v_n(a))^{\top}$ satisfy the Jordan-Pochhammer
systems
\begin{eqnarray}\label{JP2}
du=\Omega^u\, u, \qquad dv=\Omega^v\, v,
\end{eqnarray}
where
\begin{eqnarray*}
\Omega^u=\sum_{1\leqslant j<k\leqslant n}J_{jk}(\beta^u)\,
\frac{d(a_j-a_k)}{a_j-a_k}, & \qquad \beta_i^u=\beta_i^1-\beta_i^2, \\
\Omega^v=\sum_{1\leqslant j<k\leqslant n}J_{jk}(\beta^v)\,
\frac{d(a_j-a_k)}{a_j-a_k}, & \qquad
\beta_i^v=\beta_i^2-\beta_i^3.
\end{eqnarray*}
Therefore for the functions $u_i(a)$, $v_i(a)$ there hold integral
representations similar to (\ref{JP}). At the same time, the
functions $b_i(a)$ satisfy the following system of non-homogeneous
linear differential equations:
\begin{eqnarray}\label{schl3}
db_i(a)=-\sum_{j=1,j\ne i}^n(\beta_ib_j(a)-\beta_jb_i(a))\,
\frac{d(a_i-a_j)}{a_i-a_j}+\theta_i, \qquad i=1,\ldots,n,
\end{eqnarray}
where
$$
\theta_i=-\sum_{j=1,j\ne i}^n(u_i(a)v_j(a)-u_j(a)v_i(a))\,
\frac{d(a_i-a_j)}{a_i-a_j}, \qquad \beta_i=\beta_i^1-\beta_i^3.
$$

Denote by $b(a)$ the vector function
$(b_1(a),\ldots,b_n(a))^{\top}$ and by $\theta$ the vector valued
differential 1-form $(\theta_1,\ldots,\theta_n)^{\top}$, and
rewrite the system (\ref{schl3}) in the matrix form
\begin{eqnarray}\label{JP3}
db=\Omega\,b+\theta,
\end{eqnarray}
where the matrix 1-form $\Omega$ has the form (\ref{JPform}). For
any solutions $u(a)$, $v(a)$ of the systems (\ref{JP2}), the
system (\ref{JP3}) is integrable in the polydisk $D(a^0)$, which
follows from the integrability of the Schlesinger equation. Hence
its solution can be represented in the form $b(a)=Y(a)c(a)$, where
$Y(a)$ is a fundamental matrix of the corresponding homogeneous
system and the vector function $c(a)$ satisfies the relation
$$
dc=Y^{-1}\theta,
$$
whose right part is necessarily a closed differential 1-form.
Therefore,
$$
c(a)={\rm const}+\int_{a^0}^aY^{-1}\theta.
$$
and the integral does not depend on a path connecting the points
$a^0$ and $a\in D(a^0)$.

In the (generic) case, when all $\beta_i\ne 0$, for the entries of
the matrix $Y^{-1}$ one can propose integral representations
analogous to (\ref{JP}). Namely, the following lemma holds.
\medskip

{\bf Lemma 2} (A.\,Varchenko). {\it Let all $\beta_i\ne0$ and
$Y^*(a)$ be a fundamental matrix of a linear Pfaffian system
\begin{eqnarray}\label{JPt}
db^*=-\Omega\,b^*,\qquad -\Omega=\sum_{1\leqslant j<k\leqslant
n}J_{jk}(-\beta)\,\frac{d(a_j-a_k)}{a_j-a_k}.
\end{eqnarray}
Then $Y^{-1}=(Y^*)^{\top}\Lambda$, where $\Lambda={\rm diag}\,
(1/\beta_1,\ldots,1/\beta_n)$.}
\medskip

{\bf Proof.} Since for a pair of vectors $u,w\in{\mathbb C}^n$ the
scalar products $\Bigl(u,J_{jk}(\beta)w\Bigr)$ and
$\Bigl(J_{jk}(\beta)u,w\Bigr)$ are equal to
\begin{eqnarray*}
\Bigl(u,J_{jk}(\beta)w\Bigr)&=&u_j(\beta_kw_j-\beta_jw_k)+u_k(\beta_jw_k-\beta_kw_j),\\
\Bigl(J_{jk}(\beta)u,w\Bigr)&=&w_j(\beta_ku_j-\beta_ju_k)+w_k(\beta_ju_k-\beta_ku_j),
\end{eqnarray*}
then for the symmetric bilinear form
$$
(u,w)_{\beta}=\sum_{i=1}^n\frac1{\beta_i}\,u_iw_i
$$
one has $\Bigl(u, J_{jk}(\beta)w\Bigr)_{\beta}=
\Bigl(J_{jk}(\beta)u, w \Bigr)_{\beta}$ for all $j<k$, or
$$
\Bigl(J_{jk}(\beta)u, w \Bigr)_{\beta}+\Bigl(u, J_{jk}(-\beta)w
\Bigr)_{\beta}=0.
$$
Consequently, if $b(a)$ is a solution of the system
(\ref{JPsyst}), (\ref{JPform}), and $b^*(a)$ is a solution of the
system (\ref{JPt}), then the function $\Bigl(b(a), b^*(a)
\Bigr)_{\beta}$ is constant. Therefore,
$$
(Y^*)^\top\Lambda\,Y=C
$$
is a constant non-degenerate matrix (choosing a suitable
fundamental matrix $Y$, one can suppose that the matrix $C$ is the
identity matrix). {\hfill $\Box$}

\bigskip

Institute for Information Transmission Problems, Russian Academy of Sciences,
19, Bolshoy Karetniy, 127994, Moscow, Russia; rgontsov@inbox.ru
\\

Moscow State Regional Social-Humanitarian Institute, 30, Zelenaya, 
140411, Kolomna, Russia; lexin\_vp@mail.ru

\end{document}